\documentclass[review]{elsarticle}
\RequirePackage{amssymb}
\usepackage{bm}
\RequirePackage{amsmath}
\usepackage{amsthm}
\usepackage{color}
\usepackage{lineno,hyperref}
\modulolinenumbers[5]
\usepackage{geometry}
\usepackage{lineno,hyperref}
\modulolinenumbers[5]

\journal{Journal of Graph Theory}

\begin{document}
\newtheorem{theorem}{Theorem}[section]
\newtheorem{definition}{Definition}[section]
\newtheorem{lemma}{Lemma}[section]
\newtheorem{claim}{Claim}[section]
\newtheorem{corollary}{Corollary}[section]
\newtheorem{proposition}{Propotion}[section]
\newtheorem{conjecture}{Conjecture}[section]
 \abovedisplayskip 6pt plus 2pt minus 2pt \belowdisplayskip 6pt
 plus 2pt minus 2pt
 %%%%%%%%%%%%%%%%
\def\vsp{\vspace{1mm}}
\def\th#1{\vspace{1mm}\noindent{\bf #1}\quad}
\def\proof{\vspace{1mm}\noindent{\it Proof}\quad}
\def\no{\nonumber}
\newenvironment{prof}[1][Proof]{\noindent\textit{#1}\quad }
{\hfill $\Box$\vspace{0.7mm}}
\begin{frontmatter}

\title{An Improved Error Term for Tur$\acute{\rm a}$n Number of Expanded Non-degenerate 2-graphs}
\author[nuaa]{Yucong Tang}
\ead{tangyucong@amss.ac.cn}

\author[niu]{Xin Xu\fnref{NSF2}}
\ead{xuxin@ncut.edu.cn}
\fntext[NSF2]{Research partly supported by Natural Science Foundation of Beijing (Grant No. 1174015)}

\author[ams,ucas]{Guiying Yan\corref{cor1}\fnref{NSF1}}
\ead{yangy@amss.ac.cn}
\fntext[NSF1]{Research partly supported by  National Natural Science Foundation of China (Grant No. 11631014)}

\cortext[cor1]{Corresponding author}
\address[ams]{Academy of Mathematics and Systems Science, Chinese Academy of
Sciences,
 \\
 Beijing $100190$, P. R. China}

\address[niu]
{School of Sciences, North China University of Technology\\ Beijing $100144$, P. R. China}

\address[nuaa]
{College of Science, Nanjing University of Aeronautics and Astronautics, Nanjing $211106$, P. R. China}

\address[ucas]
{School of Mathematical Sciences, University of Chinese Academy of Sciences, Beijing $100049$, P. R. China}

\begin{abstract}
 For a 2-graph $F$, let $H_F^{(r)}$ be the $r$-graph obtained from $F$ by enlarging each edge with a new set of $r-2$ vertices.  We show that if $\chi(F)=\ell>r \geq 2$, then $ {\rm ex}(n,H_F^{(r)})= t_r (n,\ell-1)+ \Theta( {\rm biex}(n,F)n^{r-2}),$ where  $t_r (n,\ell-1)$ is the number of edges of an $n$-vertex complete balanced $\ell-1$ partite $r$-graph and ${\rm biex}(n,F)$ is the extremal number of the decomposition family of $F$.  Since ${\rm biex}(n,F)=O(n^{2-\gamma})$ for some $\gamma>0$, this improves on the bound ${\rm ex}(n,H_F^{(r)})= t_r (n,\ell-1)+ o(n^r)$ by Mubayi (2016) \cite{Mubayi&Dhruv2016}. Furthermore, our result implies that  ${\rm ex}(n,H_F^{(r)})= t_r (n,\ell-1)$ when $F$ is  edge-critical, which is an extension of the result of Pikhurko (2013) \cite{Pikhurko}.
\end{abstract}


\end{frontmatter}

\section{Introduction}
An $r$-graph (or  $r$-uniform hypergraph) $G$ consists of a vertex set and an edge set with exactly $r$ vertices in each edge. We sometimes identify an $r$-graph $H$ with its edge set, and denote by $V(H)$
 its vertex set. An $r$-\emph{clique} of \emph{order} $j$, denoted by $K_j^{(r)}$, is an $r$-graph on $j\geq r$ vertices consisting of all ${j\choose r}$ different $r$-tuples. An $r$-graph $G$ is said to be $k$-partite if its vertex set can be partitioned into $k$ classes $V_1\cup V_2\cup \cdots\cup V_k$ such that every edge of $G$ contains at most one vertex in $V_i$, $i=1,\ldots,k$. We say $G$ is a complete $k$-partite $r$-graph if $G$ consists of all $r$-tuples  intersecting each vertex class in at most 1 vertex. Given two $r$-graphs $G$ and $F$, we say $G$ is $F$-free if $G$ does not have a (not necessarily induced) subgraph isomorphic to $F$. For a positive integer $\ell$, we denote by $[\ell]$ the set $\{1,\ldots,\ell\}$. For a set $V$ and an integer $r\geq 1$, let $[V]^r$ be the set of all $r$-element subsets of $V$. We write $[\ell]^r$ instead of $[[\ell]]^r$ for simple.

 The Tur$\acute{\rm a}$n number ${\rm ex}(n,F)$ is the maximum number of edges in an $n$-vertex $F$-free $r$-graph.
A simple and important averaging argument of Katona, Nemetz and Simonovits \cite{Katona} shows that ${n\choose r}^{-1}{\rm ex}(n,F)$ form a decreasing sequence of real numbers in $[0,1]$. It follows that the sequence has a limit, called the Tur$\acute{\rm a}$n density and denoted by $\pi(F)$.

The Tur$\acute{\rm a}$n density is an asymptotic result ${\rm ex}(n,F)\sim \pi(F){n\choose r}$. An important fundamental theorem proved by  Erd$\ddot{\rm o}$s and Stone \cite{Erdos} characterizes the Tur$\acute{\rm a}$n density of any 2-graph with its chromatic number.
 \begin{theorem}[{\cite{Erdos}}]\label{ES}
 Let $F$ be a 2-graph with $\chi(F)=\ell$, then $\pi(F)=\pi(K_{\ell})=\frac{\ell-2}{\ell-1}$.
 \end{theorem}

However, when $F$ is an $r$-graph, $\pi(F)\neq 0$, and $r>2$, determining $\pi(F)$ is a hard problem, even for very simple $r$-graphs. In this paper, we focus on the non-degenerated expanded 2-graphs which are $r$-graphs defined as follows. Let $\ell> r\geq 2$ and let $F$ be a 2-graph, and $H_F^{(r)}$ be the $r$-graph obtained from $F$ by enlarging each edge with a new set of $r-2$ vertices.  Mubayi \cite{Mubayi1} first determined the Tur$\acute{\rm a}$n density of expanded cliques and obtained a stability result.

\begin{theorem}[{\cite{Mubayi1}}]\label{th2}
 Let $\ell> r\geq 2$. Then
 \begin{align*}
 \pi(H_{K_{\ell}}^{(r)})=\frac{r!}{(\ell-1)^r}{{\ell-1}\choose r}.
 \end{align*}
\end{theorem}
Later, using the stability method, Pikhurko \cite{Pikhurko} obtained the exact number of ${\rm ex}(n, H_{K_{\ell}}^{(r)})$.

It was mentioned in the survey of Mubayi \cite{Mubayi&Dhruv2016} that   Alon and Pikhurko  observed that the approach applied to prove Theorem 2 in \cite{Pikhurko} can
be extended to any edge-critical graph $F$ with $\chi(F) > r$. More generally, the following results can be easily achieved through a result of Erd$\ddot{\rm o}$s \cite{Erdos1964},
the supersaturation technique (see Erd$\ddot{\rm o}$s-Simonovits \cite{Erdos&Simonovits1983}), and Theorem \ref{th2}.

\begin{theorem}\label{Main}
Let $ \ell> r\geq 2 $. Let  $F$  be any 2-graph with $\chi(F)=\ell$, then
 \begin{align*}
 \pi(H_{F}^{(r)})=\frac{r!}{(\ell-1)^r}{{\ell-1}\choose r}.
 \end{align*}
\end{theorem}

This asymptotically gave  the Tur$\acute{\rm a}$n number of all non-degenerated expanded 2-graphs.

\vspace{2mm}For self completeness, we will give a short proof of Theorem \ref{Main} in the next section.
Then we prove a stability result of $H_{F}^{(r)}$.
Denote by $T_r(n,\ell)$ the complete $\ell$-partite $r$-graph on $n$ vertices, where the size of   each vertex class differs at most 1, and  set $t_r(n,\ell)=|T_r(n,\ell)|$. We say two $r$-graphs $G_1$ and $G_2$ of  order $n$ are $\varepsilon$-close if we can add or remove at most $\varepsilon {n \choose r}$ edges
from   $G_1$ to make it isomorphic to $G_2$; in other words, for some bijection
$\sigma: V (G_1)\rightarrow V(G_2)$ the symmetric difference between $\sigma(G_1) = \{\sigma(D): D \in G_1 \}$ and $G_2$ has at most $\varepsilon {n \choose r}$ edges.

\begin{theorem}[Stability of $H_F^{(r)}$]\label{Stability} Fix $\ell > r\geq 2$ and 2-graph $F$ with $\chi(F)=\ell$. For every $\varepsilon>0$, there exist $n_0 = n_0(r, \ell, \varepsilon)>0$, $\eta = \eta(r, \ell, \varepsilon) > 0$, such that if $n>n_0$ and  $G$ is an $n$-vertex $H_F^{(r)}$-free $r$-graph  with  $|G|\geq |T_r (n,\ell-1)|-\eta {n\choose r}$, then $G$ is $\varepsilon$-close to $T_r(n,\ell-1)$.
\end{theorem}

\begin{definition}[${\rm biex}(n,F)$]
Given a 2-graph $F$ with $\chi(F)=\ell$, the decomposition family $\mathcal{F}_F$ of  $F$ is the set of bipartite graphs which are obtained from $F$ by deleting $\ell-2$ colour classes in some $\ell$-colouring of $F$. Observe that $\mathcal{F}_F$ may contain graphs which are disconnected, or even have isolated vertices. Let $\mathcal{F}_F^{\star}$ be a minimal subfamily of $\mathcal{F}_F$ such that for any $H\in \mathcal{F}_F$, there exists $H' \in \mathcal{F}_F^{\star}$ with $H' \subset H$. We define
$${\rm biex}(n,F):={\rm ex}(n,\mathcal{F}_F)={\rm ex}(n,\mathcal{F}_F^{\star}).$$
%Further more, we define $${\rm biex}(n,H_F^{(r)}):={\rm ex}(n,\mathcal{F}_F^{(r)})={\rm ex}(n,\mathcal{F}_F^{\star (r)}).$$
%where $\mathcal{F}_F^{(r)}$ is the expanded family of $\mathcal{F}_F$, i.e., $\mathcal{F}_F^{(r)}=\{H_H^{(r)}:H\in\mathcal{F}_F\}$.
\end{definition}

Furthermore, we prove an improved bound for the Tur$\acute{\rm a}$n function ${\rm ex}(n,H_F^{(r)})$.

\begin{theorem}\label{turannum}
%Given $2$-graph $F$ with $\chi(F)=\ell>r$, there exist $c_1,c_2,n_0>0$ such that if $n\geq n_0$, we have
%$$t_r (n,\ell-1)+c_1 {\rm biex}(n,H_F^{(r)})\leq {\rm ex}(n,H_F^{(r)})\leq t_r (n,\ell-1)+c_2 {\rm biex}(n,F)n^{r-2}.$$
Given a $2$-graph $F$ with $\chi(F)=\ell>r$, then
$$ {\rm ex}(n,H_F^{(r)})= t_r (n,\ell-1)+ \Theta( {\rm biex}(n,F)n^{r-2}).$$
\end{theorem}

We say a 2-graph $F$ is \emph{edge-critical} if there exists an edge $e \in E(F)$ such that $\chi(F) > \chi(F-e)$. The following theorem is a direct corollary of Theorem \ref{turannum}.

\begin{theorem}\label{edgecritical}
Given a $2$-graph $F$ with $\chi(F)=\ell>r$. If $F$ is edge-critical, then
$$  {\rm ex}(n,H_F^{(r)})= t_r (n,\ell-1).$$
\end{theorem}

In the following section, we will a short proofs  of Theorem \ref{Main} by using hypergraph Lagrange method. In section 3, we prove Theorem \ref{Stability} based on  the hypergraph removal lemma and a stability result of expanded cliques (see \cite{Pikhurko}). In the last section, we prove Theorem \ref{turannum}, the idea  is  first to identify a  copy of $F$ in the $'$2-shadow$'$ of an $r$-graph and then extend this copy to $H_F^{(r)}$.

\section{Proof of Theorem \ref{Main}}

The hypergraph Lagrange method  was developed independently by Sidorenko \cite{SidorenkoL} and Frankl and F$\ddot{\rm u}$redi \cite{FFL}.

Let $G$ be an $r$-graph on $[n]=\{1,\ldots, n\}$ with edge set $E$, and $x=\{x_1,\ldots,x_n\}$ with $x_i\geq 0$ for all $1\leq i\leq n$ and $\sum_{i=1}^n x_i=1$. Define
$$p_G(x)=\sum_{\{i1,i2,\ldots,ir\}\in E}x_{i1}x_{i2}\cdots x_{ir},$$

The Lagrange of $G$ is defined as $\lambda(G)=\max_x p_G(x)  $. $G$ is said to be dense if the inequality $\lambda(G')< \lambda(G)$ holds for all its proper subgraphs $G'\neq G$.
We say that $G$   \emph{covers pairs} if for every pair of vertices $i,j$ in $G$, there is an edge containing both $i$ and $j$. We need the following results.

\begin{lemma}[{\cite{SidorenkoL},\cite{FFL}}]\label{dense}
Every dense graph covers pairs.
\end{lemma}

Given $r$-graphs $F$ and $G$, we say $f:V(F)\rightarrow V(G)$ is a \emph{homomorphism} if $f(e)\in E(G)$ for all $e\in E(F)$. And we call $G$ is $F$-hom-free if there is no homomorphism from $F$ to $G$. The following theorem shows how to compute the Tur$\acute{\rm a}$n density of any $r$-graph.

\begin{lemma}[{\cite{SidorenkoL}}]\label{HL}
Let $F$ be an $r$-graph, then
 $$\pi(F)=\sup\{r! \lambda(G):{\mbox{$G$ is $F$-hom-free}}\}.$$
\end{lemma}

\begin{prof}[Proof 1]
First, for the $(\ell-1)$-clique  $K_{\ell-1}^{(r)}$, we have that $K_{\ell-1}^{(r)}$ is $H_F^{(r)}$-hom-free since $\chi(F)=\ell$. Otherwise, there  is a homomorphism $f: F\rightarrow K_{\ell-1}^{(r)}$,  then $f^{-1}(v_i)$ form  a vertex partition of $F$ and every edge of $H_F^{(r)}$ contains at most one vertex in $f^{-1}(v_i)$. Thus $f^{-1}(v_i)$ is an independent set in $F$, which is a contradiction.

On the other hand, for any dense graph $G$ on at least $\ell$ vertices, we can construct a homomorphism $g$ from   $H_F^{(r)}$ to $G$. Since $\chi(F)=\ell$, so there is a partition of $V(F)$ into $\ell$ independent set. We map   each independent set  to   $\ell$ distinct vertices of $G$. For the rest vertices in $H_F^{(r)}$, we denote the vertices in the edge containing $i,j$ by $v_{ij}^k$, $k=1,\ldots,r-2$. By Lemma \ref{dense}, $G$ covers pairs. So there is an edge containing both $g(i)$ and $g(j)$.  We  map  $v_{ij}^k$ to the rest $r-2$ distinct vertices in that edge. Thus $g$ is a homomorphism and by Lemma \ref{HL}, we have
$$\pi(F)=\sup_{\mbox{$G$ is $F$-hom-free}}r! \lambda(G)= r!\lambda(K_{\ell-1}^r)=\frac{r!}{(\ell-1)^r}{(\ell-1)\choose r},$$
which complete the proof.
\end{prof}

\section{Stability of $H_F^{(r)}$}

To proof Theorem \ref{Stability},
we   need the following stability result of $H_{K_{\ell}}^{(r)}$ and the hypergraph removal lemma.

\begin{lemma}[{\cite{Pikhurko}}]\label{HKLStab}
Fix $\ell > r \geq 2$. For every $\varepsilon_1  > 0$, there are $\eta_1 = \eta_1(r, \ell, \varepsilon_1) > 0$
and $n_0 = n_0(r, \ell, \varepsilon_1)$ such that any $H_{K_{\ell}}^{(r)}$-free $r$-graph $G$ of order $n \geq n_0$ and size at least $|T_r (n,\ell-1)|-\eta_1 n^r$
is
$\varepsilon_1$-close to $T_r(n, \ell-1)$.
\end{lemma}

Hypergraph removal lemma is yield  among a series extensions of the Szemer$\acute{\rm e}$di's regularity lemma to $r$-graphs (see \cite{10,20,27,Regularity}). Tao \cite{Tao}   also obtained such a generalization.    In this paper, we will use two versions of the hypergraph removal lemma as follows.
\begin{lemma}[Hypergraph Removal Lemma, \cite{Regularity}]\label{RemovalI}
Fix an $r$-graph $F$. For every $\varepsilon>0$, there exist $\eta>0$ and $n_0>0$ such that for every $n$-vertex  $r$-graph $G$ with $n>n_0$, if $G$ contains at most $\eta n^{|V(F)|}$ copies of  $F$, then  one can delete at most $\varepsilon {n\choose r}$ edges to make it $F$-free.
\end{lemma}

The second version is as follows. The proof is also based on hypergraph regularity lemma and general dense counting lemma and similar to that of Lemma \ref{RemovalI}.
\begin{lemma}\label{Removal}
Fix an $r$-graph $F$. For every $\eta_2>0$, there exist $n_0>0$ such that for every $n$-vertex  $r$-graph $G$ with $n>n_0$, if $G$ is $F$-free, then  one can delete at most $\eta_2 {n\choose r}$ edges to make it $F$-hom-free.
\end{lemma}

\begin{prof}[Proof of Theorem \ref{Stability}]We choose constant $\varepsilon_1+\eta_2<\varepsilon$ and  $\eta+\eta_2<\eta_1$.

According to Lemma \ref{Removal}, we can delete at most $\eta_2 {n\choose r}$ edges, denoting the remain $r$-graph by $G'$, to make $G'$ $H_F^{(r)}$-hom-free, which implies  $G'$ is $H_{K_{\ell}}^{(r)}$-free, and
$$|G'|\geq |G|-\eta_2 {n\choose r}\geq |T_r (n,\ell-1)|-\eta {n\choose r}-\eta_2 {n\choose r}\geq |T_r (n,\ell-1)|-\eta_1{n\choose r}. $$

Apply Lemma \ref{HKLStab} to $G'$ for  $\varepsilon_1$, we have $G'$ is $\varepsilon_1$-close to $T_r (n,\ell-1)$. Thus $G$ is $(\varepsilon_1+\eta_2)$-close to $T_r (n,\ell-1)$, which complete the proof.
\end{prof}

\section{Proof of Theorem \ref{turannum}}
For real constants $\alpha,\beta$, and a non-negative constant $\xi$, we write $$\alpha=\beta\pm \xi, \quad \mbox{if\ } \beta-\xi\leq \alpha \leq \beta+\xi.$$

 For $U\subseteq V(H^{(r)})$, we denote by $H^{(r)}[U]$ the sub-hypergraph of $H^{(r)}$ induced on $U$ (i.e. $H^{(r)}[U]=H^{(r)}\cap [U]^r$).

 Given vertex sets $V_1,\ldots, V_{\ell}$, let $K_{\ell}^{(j)}(V_1,\ldots,V_{\ell})$ be the complete $\ell$-partite, $j$-graph. If $|V_i|=m$ for all $i\in [\ell]$, then an $(m,\ell,j)$-graph $H^{(j)}$ on $V_1\cup\cdots\cup V_{\ell}$ is any subset of $K_{\ell}^{(j)}(V_1,\ldots,V_{\ell})$. Also, we regard the vertex partition $V_1\cup\cdots\cup V_{\ell}$ as an $(m,\ell,1)$-graph $H^{(1)}$. For $j\leq i\leq \ell$ and set $\Lambda_i\in [\ell]^i$, we denote the $\cup_{\lambda\in \Lambda_i}V_{\lambda}$ induced sub-hypergraph of the $(m,\ell,j)$-graph $H^{(j)}$ by $H^{(j)}[\Lambda_i]=H^{(j)}[\cup_{\lambda\in \Lambda_i}V_{\lambda}]$.

To prove Theorem \ref{turannum}, it is sufficient to prove the following theorem.

\begin{theorem}\label{4.1}
Given a $2$-graph $F$ with $\chi(F)=\ell>r\geq 2$, there exist $c_1,c_2,n_0>0$ such that if $n\geq n_0$, we have
$$t_r (n,\ell-1)+c_1 {\rm biex}(n,F)n^{r-2}\leq {\rm ex}(n,H_F^{(r)})\leq t_r (n,\ell-1)+c_2 {\rm biex}(n,F)n^{r-2}.$$
\end{theorem}

\vspace{1mm}\begin{prof}[Proof of Theorem \ref{4.1}]
Firstly, the left hand-side inequality  is obtained as follows. Let $H$ be an $n$-vertex $\mathcal{F}_F$-free $2$-graph with $ {\rm biex}(n,F)$ edges, and let $c=(\ell-1)^{-2}$, $c_1={{\ell-2}\choose {r-2}}\frac{c }{(\ell-1)^{r-2}}$. Obviously, there exists an $n/(\ell-1)$-vertex subgraph $H'$ of $H$ with at least $c¡¤|H|$ edges.

 Next, we construct $G$  from $T_r(n,\ell-1)$ as follows. Without loss of generality, let $V_1$ be the vertex class  of $T_r(n,\ell-1)$ with largest size, we insert $H'$ into $V_1$. Then for each edge $(u,v)$ in $H'$,  add all the $r$-tuples that contains $u,v$ and  $(r-2)$ vertices chosen from different vertex classes except $V_1$ to $G$, i.e., $$G=T_r(n,\ell-1)\cup (\cup_{(u,v)\in H'}E(u,v)) ,$$
where $E(u,v)=\{\{u,v\}\cup f: |f|=r-2,|f\cap V_i|\leq 1, f\cap V_1=\emptyset,i=2\ldots,\ell-1\}$.

Clearly, we have
$$|G|\geq t_r(n,\ell-1)+c{{\ell-2}\choose {r-2}}\frac{1}{(\ell-1)^{r-2}} {\rm biex}(n ,F)n^{r-2},$$ and by definition of $\mathcal{F}_F$, the graph $G$ is $H_F^{(r)}$-free, and therefore $$t_r (n,\ell-1)+c_1 {\rm biex}(n ,F)n^{r-2}\leq {\rm ex}(n,H_F^{(r)}).$$

 Secondly,  the main idea to prove the right hand-side inequality is to find a copy of $F$ in the 2-shadow of $G$, and then extend $F$ to $H_F^{(r)}$ in $G$. Here by saying 2-shadow of $G$, denoted by $\Delta_2(G)$, we mean the set of all 2-tuples $\{u,v\}\in [V(G)]^2$ that are contained in some edge of $G$.

 Set $|V(F)|=m$ and choose $\varepsilon>0$ small enough. Suppose $G$ is an $n$-vertex $H_F^{(r)}$-free $r$-graph with $|G|> t_r (n,\ell-1)+c_2 {\rm biex}(n,F)n^{r-2}$, $n\geq n_0$, then by Theorem \ref{Stability} $G$ is $\varepsilon$-close to $T_r (n,\ell-1)$. Thus  $V(G)$  can be  partitioned  into balanced $V_1\dot{\cup} \cdots\dot{\cup} V_{\ell -1}$ corresponding to $T_r(n,\ell-1)$.

Since $\frac{n}{2}\leq{\rm biex}(n,F)=o(n^2)$ or ${\rm biex}(n,F)=0$, so we have

\vspace{2mm}\noindent {\bf Fact 1.}\quad $\frac{1}{2}n^{r-1}<{\rm biex}(n,F)n^{r-2}=o(n^r)$ or ${\rm biex}(n,F)n^{r-2}=0$.

\vspace{2mm}

We call a pair of vertices \emph{bad} if it is covered by at most
 $$\kappa(n,F)=|V(H_F^{(r)})|{n\choose {r-3}}$$
edges of $G$.

Let $G'$ be obtained from $G$ by deleting  all edges containing bad pairs, at most ${n\choose 2}\kappa(n,F)<\varepsilon{n\choose r}$. So $G'$ is   $2\varepsilon$-close to  $T_r (n,\ell-1)$.

For any vertex $v$, we denote by $d(v)$ the vertex degree of $v$ in $G'$, and  denote by $N(v)$  the neighbours of $v$ in $G'$, i.e., for each vertex $u$ in $N(v)$, there is an  edge containing both $u$ and  $v$. Let $N_{V_i}(v)=N(v)\cap V_i$ and $d_i(v)=|N_{V_i}(v)|$. An edge $e$ is \emph{crossing} if $|e\cap V_i|\leq 1$ for $i\in [\ell-1]$. Let $C(v)$ be the set of crossing edges   containing $v$ and $\bar{C}(v)$ be the set of non-crossing edges containing $v$, and we call $d^c(v)=|C(v)|$ the \emph{crossing degree} of $v$ and   $\bar{d}^c(v)=|\bar{C}(v)|$ the non-crossing degree.

Observe first that we may assume without loss of generality that
   \begin{align}\label{1}\delta (G') \geq \delta(T_{r}(n,\ell-1)).\end{align}
where $\delta (G')=\min\{d(v):v\in V(G')\}$.
Indeed, if this is not the case, we can repeatedly delete vertices of minimum degree  of $G'$ and delete all edges containing bad pairs until we arrive at a graph $G'_{n^{\star}}$ on $n^{\star}$ vertices with $\delta (G_{n^{\star}}) \geq \delta(T_{r}(n^{\star},\ell-1))$. Denote the sequence of graphs obtained in this way by $G'_n:=G',G'_{n-1},\ldots ,G'_{n^{\star}}$. We need to verify that $n^{\star}\geq n_0$.  Indeed, we have
\begin{align*}|G'_{n-1}|&\geq |G'_n|-\delta(G'_n)-{n\choose 2}\kappa(n,F)\\
 &> |T_r(n,\ell-1)|-\delta(T_r(n,\ell-1))+\bigg(c_2{\rm biex}(n,F)n^{r-2}-{n\choose 2}\kappa(n,F)\bigg)\\
 &\geq |T_r(n-1,\ell-1)|+\frac{c_2}{2} {\rm biex}(n-1,F)(n-1)^{r-2}+1.\end{align*}
Similarly, we have
$$|G'_{n-i}|\geq |T_r(n-i,\ell-1)|+\frac{c_2}{2^i} {\rm biex}(n-i,F)(n-i)^{r-2}+i.$$
Let $i^{\star}:= n-n^{\star}$. If $n^{\star}<n_0$, then $i^{\star}>n-n^{\frac{1}{2r}}$, which implies $i^{\star}\geq {{n-i^{\star}}\choose r}$, a contradiction. Hence we may assume (\ref{1}).

Next we move the vertices   to get a max $(\ell-1)$-cut of $G'$, i.e., maximise the number of crossing edges. For $v\in V_i$ and $i\neq j\in [\ell-1]$, let
$$E_{s,t}^{i,j}(v)=\{e\in C(v)\cup \bar{C}(v): |e\cap V_i|=s, |e\cap V_j|=t, |e\cap V_k|\leq 1, k\neq i,j\}$$
Then, the max $(\ell-1)$-cut implies  a vertex partition such that for each vertex $v\in V_i$, we have \begin{align}\label{4}
|E_{2,0}^{i,j}(v)|\leq |E_{1,1}^{i,j}(v)| \quad j\neq i, j\in [\ell-1]
\end{align}

Since the number of crossing edges is at least $t_r (n,\ell-1)-2\varepsilon{n\choose r}$, so a simple computation would indicate that \begin{align}    |V_i|=(1\pm \sqrt{\varepsilon})\frac{n}{\ell-1},\quad i\in [\ell-1]\label{5}\end{align}

Note that $d^c(v)+\bar{d}^c(v)=d(v)$. Let $X=\{v\in V(G'): \bar{d}^c(v)> \varepsilon^{1/(4(r-1))}n^{r-1}\}$ and $X_i=X\cap V_i$.
Set $V'_i=V_i\backslash X$. Since $G'$ is $2\varepsilon$-close to $T_r(n,\ell-1)$, so $$ \frac{1}{r}|X|\varepsilon^{1/(4(r-1))}n^{r-1}\leq 2\varepsilon {n\choose r},$$
which implies $|X|\leq \varepsilon^{2/3}n$.

Then for every $v\in V'_i$, $i\in [\ell-1]$, we have
\begin{align} \label{8}d^c(v)&=d(v)-\bar{d}^c(v)\geq \delta(T_r(n,\ell-1))- \varepsilon^{1/(4(r-1))}n^{r-1}
\end{align}
and this implies that for every $j\neq i,j\in [\ell-1]$,
\begin{align}
|N_{V_j}(v)\backslash X| &\geq |V_i|-|X|-(\ell-1)^{r-2}\varepsilon^{\frac{1}{(4(r-1))}}n \notag \\
&\geq(1-2(\ell-1)^{r-1}\varepsilon^{\frac{1}{(4(r-1))}})\frac{n}{\ell-1}\notag\\
&\geq(1-\varepsilon^{\frac{1}{(5(r-1))}})\frac{n}{\ell-1}\label{2}
\end{align}

Let $q$ be a positive constant depending only on $|V(F)|$ and $\varepsilon $. Its value will be given later.

\vspace{2mm}\noindent{\bf Case 1.} If $|X|< q(m-1)$ and ${\rm biex}(n,F)n^{r-2}>0$.
Since $|G|> t_r (n,\ell-1)+c_2 {\rm biex}(n,F)n^{r-2}$, so the number of non-crossing edges in $G'$ is at least \begin{align}\label{crossingedgenum}\frac{c_2}{2} {\rm biex}(n,F)n^{r-2}.
\end{align} For every $i\in [\ell-1]$, we denote by $E(V'_i)$ the set of non-crossing edges in $G'$ that contains at least 2 vertices in $V'_i$. We have
$$\sum_{i=1}^{\ell-1}|E(V'_i)|\geq \frac{c_2}{2} {\rm biex}(n,F)n^{r-2}-|X|n^{r-1}\geq \frac{c_2}{3}{\rm biex}(n,F)n^{r-2},$$
where the last inequality is due to Fact 1 and $n$ sufficiently large. Then, there exists some $i^{\star}$ such that \begin{align}|E(V'_{i^{\star}})|\geq {\rm biex}(n,F)n^{r-2}.\label{3}\end{align}
Next we write $D=\{\{u,v\}\in [V_{i^{\star}}]^2: \mbox{there exists some } e\in E(V_{i^{\star}}),\ \mbox{such that } \{u,v\}\subseteq e\}$. Because each vertex pair in $D$ is contained in at most $n^{r-2}$ edges in $E(V_{i^{\star}})$, so, by (\ref{3}), we have
$|D|\geq {\rm biex}(n,F)$. That is, we can find some $H\in \mathcal{F}^{\star}_F$ in $\Delta_2(G')[V'_{i^{\star}}]$. Let such a copy of $H$ be fixed and assume without loss of generality that $V(H)\subseteq V_{\ell-1}$. Then we show that $H$ can be extended to a copy of $F$ in the 2-shadow of $G'$ by finding a complete $(\ell-2)$-partite 2-graph in $\Delta_2(G')$.

Note that by (\ref{2}), we have for any vertex set $S\subseteq V(G)$ with $S\leq \ell m$ and every $i\in[\ell-1]$, the number of common neighbours in $V'_i$ of every vertex in $S$ is at least
\begin{align}
(1-\ell m\varepsilon^{\frac{1}{(5(r-1))}})\frac{n}{\ell-1}-\ell m-|X|\geq(1- \varepsilon^{\frac{1}{(6(r-1))}})\frac{n}{\ell-1}.
\end{align}
The   inequality is due to $\varepsilon$ is small enough and $n$ is sufficiently large.

We inductively find sets $S_i\subseteq V'_i$ of size $m$ which form the parts of the complete $(\ell-2)$-partite 2-graph. For each $1\leq i\leq \ell-2$ in turn, we note that $|V(H)|+(i-1) m\leq \ell |V(H)|$, and therefore the set $V(H)\cup S_1\cup\cdots\cup S_{i-1}$ has at least $(1- \varepsilon^{\frac{1}{6(r-1)}})\frac{n}{\ell-1}\geq m $ common neighbours in $V'_{i}$. We let $S_i$ be any set of size $m$ of these common neighbours. Hence we can extend $H$ to a copy of $F$ in $\Delta_2( G')$.

Finally, recalling that we have deleted the edges that contains bad pairs, each vertex pair (or edge) in this copy of $F$  is contained in at least $(\ell+(r-2){\ell \choose 2}){n\choose {r-3}}$ edges of $G'$. Thus we can choose, for each pair, one of these edges that is vertex-disjoint to the chosen ones to form a copy of $H_F^{(r)}$  in $G'$.

\vspace{2mm}\noindent{\bf Case 2.} If $|X|\geq q(m-1)$ and ${\rm biex}(n,F)n^{r-2}>0$.

Let $\bar{C}_i(x)$ be the subset of  $\bar{C}(x)$ that contains at 2 vertices in $V_i\backslash \{x\}$.

 If there is $i,j\in [\ell-1]$ (Notice that  $i=j$ is possible),   and $x\in X_i$ such that $|\bar{C}_j(x)|\geq  \sqrt{\varepsilon} n^{r-1}$.

We write $D=\{\{u,v\}\in [V_{j}]^2: \mbox{there exists some } e\in \bar{C}_j(x), \mbox{ s.t. } \{u,v\}\subseteq e\}$ and we claim $|D[V'_j]|\geq {\rm biex}(n,F)$. Because the number of 2-tuples in $D$ that contains at least 1 vertex in $X_i$ is at most $\varepsilon^{\frac{2}{3}}n^2.$ Thus the number of edges in $\bar{C}_j(x)$ that contains at least 1 vertex in $X_i$ is at most $\varepsilon^{\frac{2}{3}}n^2n^{r-3} $. Note that each 2-tuple in $D[V'_j]$ is contained in at most $n^{r-3}$ edges in $\bar{C}_j(x)$, thus
$$|D[V'_j]|\geq \frac{1}{n^{r-3}} (\sqrt{\varepsilon}-\varepsilon^{\frac{2}{3}}) n^{r-1}\geq {\rm biex}(n,F).$$
This means that we can again find a copy of some $H\in \mathcal{F}_F$. And the extending from $H$ to $F$, and then to $H_F^{(r)}$ is the same as that in Case 1.

Otherwise for every $i,j\in [\ell-1]$, $x\in X_i$, we have $|\bar{C}_j(x)|\leq  \sqrt{\varepsilon} n^{r-1}$.
Denote by $\bar{C}_i^1(x)$ the subset of $\bar{C}(x)$ that contains exactly 1 vertex in $V_i\backslash x$. Then by $\sum_{j=1}^{\ell-1} |\bar{C}_j(x)|+|\bar{C}_i^1(x)|\geq \bar{d}^c(x)$, we know
\begin{align}\label{7}|\bar{C}_i^1(x)|\geq \varepsilon^{\frac{1}{4(r-1)}} n^{r-1} -(\ell-1)\sqrt{\varepsilon} n^{r-1}\geq \varepsilon^{\frac{1}{3(r-1)}} n^{r-1}.\end{align}

\vspace{1mm}\noindent{\bf Claim 1.}\quad $|N(x)\cap V'_j|\geq \varepsilon^{\frac{1}{2(r-1)}}  n$ for every $j\in[\ell-1]$

\begin{prof}By (\ref{7}), it is easy to know that \begin{align}\label{6}|N(x)\cap V'_i|\geq \varepsilon^{\frac{1}{2(r-1)}}   n.\end{align}

Moveover, by (\ref{4}), (\ref{7}) and the assumption of $x$,  we have  for every $j\neq i$, $j\in [\ell-1]$\begin{align*}
|E_{1,1}^{i,j}(x)|\geq |E_{2,0}^{i,j}(x)|\geq |\bar{C}_i^1(x)|-(\ell-1)\sqrt{\varepsilon} n^{r-1}\geq \varepsilon^{\frac{5}{12(r-1)}} n^{r} .\end{align*}
So $|N(x)\cap V'_j|\geq \varepsilon^{\frac{1}{2(r-1)}}   n$ for $j\in [\ell-1]$.
\end{prof}

Now we   start identifying a copy of $H_F^{(r)}$ in $G'$ by 2 steps.

   Set $X':=X$. The first step is to identify $m$ vertices in $X'$ which are completely joined to an   $(r-1)$-partite $r$-graph $H_{K_{\ell-1}(m)}^{(r)}$ in $V(G')\backslash X$, with $m$  vertices, one vertex class of $K_{\ell-1}(m)$, in one vertex class of $V(G')\backslash X$. The second step is to extent the structure identified in this way  to a copy of $H_F^{(r)}$ in $G'$, which is similar as that in Case 1.

By Claim 1, we can choose for each $i$ a set $S_i \in N(x)\cap V'_i$ of size $\varepsilon^{\frac{1}{2(r-1)}}  n$. Since $G'$  is $2\varepsilon$-close to $T_r(n,\ell-1)$, so the graph $G'[\dot{\cup }S_i]$ has density at least
\begin{align*}\frac{{{\ell-1}\choose r}  (( \varepsilon^{\frac{1}{2(r-1)}}n )^{r-1}- 2\varepsilon n^{r-1})}{{{(\ell-1)\varepsilon^{\frac{1}{2(r-1)}} n}\choose r}}&\geq \frac{r!}{(\ell-1)^r}{{\ell-1}\choose r}\bigg(1-\frac{\sqrt{\varepsilon}}{2}\bigg)\\
&> \frac{r!}{(\ell-2)^r}{{\ell-2}\choose r}\\
&=\pi(H_{K_{\ell-1}}^{(r)}).\end{align*}

Then by Theorem \ref{Main}, we can not remove any edge set of size $\varepsilon n^r$ to make $r$-graph $G'[\dot{\cup }S_i]$ contain no $H_{K_{\ell-1}(m)}^{(r)}$. So by Lemma \ref{RemovalI}, there are at least $\eta n^{|V(H_{K_{\ell-1}(m)}^{(r)})|}$ copies of $H_{K_{\ell-1}(m)}^{(r)}$, where $\eta$ depends only on $m$ and $\varepsilon$. Choosing $q:=1/\eta$, we can then use the pigeonhole principle and the fact that $|X'|>q(m-1)$ to infer that there are $m$ vertices in $|X'|$ which are all adjacent to the vertices of one specific copy of $H_{K_{\ell-1}(m)}^{(r)}$ in $G'[\dot{\cup }S_i]$ as desired.

\vspace{2mm}\noindent{\bf Case 3.} When ${\rm biex}(n,F)n^{r-2}=0$, i.e., the single edge graph $H\in \mathcal{F}_F$.

The only difference is that the condition (\ref{crossingedgenum}) in Case 1 is no longer hold. We can change the proof slightly by using $G$ instead of $G'$.

 First, we change the assumption $|X| < q(m-1)$ of Case 1 to $X=\emptyset$. For $v\in V(G)$, We  call the vertex $u$ a \emph{good neighbour} of $v$ if $(u,v)$ are covered by at least $\kappa(n,F)$ edges of $G$. Note that, similar to (\ref{2}), we still have the number of good neighbours in $V_j$ of $v\in V_i$ is at least  $(1-\varepsilon^{\frac{1}{(5(r-1))}})\frac{n}{\ell-1}$. Except the only non-crossing edge we identify as a copy of $H$,  the rest of proof in Case 1 is the same.

 And in Case 2, only one vertex in $X$ is enough for us to find a copy of $H_F^{(r)}$ because there is a 1-vertex class in some coloring of $F$.

\end{prof}

\section*{References}

\bibliography{bibl}

\end{document}